\documentstyle[remark,amssymb,12pt]{article}
\newbox\smilebox
\newbox\anchorbox
\newbox\noanchorbox
\newbox\tempbox

\setbox\smilebox=\hbox{$\smile$}

\def\anchor{\hbox{\vtop{
           \hbox to \wd\smilebox{\hfil\vrule width.4pt height7pt depth1pt\hfil}
           \vskip  -11.5truept
           \hbox to \wd\smilebox{\hfil$\smile$\hfil}}}}
\setbox\anchorbox=\anchor
\def\noanchor{\hbox{\vtop{
           \hbox to \wd\anchorbox{\hfil\anchor\hfil}
           \vskip -14truept
           \hbox to \wd\anchorbox{\hfil/\hfil}}}}
\setbox\noanchorbox=\noanchor

\def\fg#1#2#3{\setbox\tempbox=\hbox{$\scriptstyle{#2}$}
\ifnum\wd\anchorbox>\wd\tempbox\dimen255=\wd\anchorbox
\else\dimen255=\wd\tempbox\fi
{#1\,\vtop{\hbox to \dimen255{\hfil\anchor\hfil}
           \vskip -6truept
           \hbox to \dimen255{\hfil$\scriptstyle{#2}$\hfil}}
           \,#3}}

\def\nfg#1#2#3{\setbox\tempbox=\hbox{$\scriptstyle{#2}$}
\ifnum\wd\noanchorbox>\wd\tempbox\dimen255=\wd\noanchorbox
\else\dimen255=\wd\tempbox\fi
{#1\,\vtop{\hbox to \dimen255{\hfil\noanchor\hfil}
           \vskip -6truept
           \hbox to \dimen255{\hfil$\scriptstyle{#2}$\hfil}}
           \,#3}}

\setbox1=\hbox{$\bot$}

\def\north#1#2{#1\,
\hbox{$\bot$\llap {\hbox to\wd1 {\hfil $/$\hfil}}}
\,#2}

\def\nao#1#2#3{#1\  \hbox{\vtop{ 
\baselineskip=4pt
\hbox{$\bot$\llap {\hbox to\wd1 {\hfil $/$\hfil}}
\hskip .05em \llap{\hbox{$^{\scriptscriptstyle{a}}$}}}\hbox{$\scriptstyle
{#2}$}}}\, #3}

\def\bp{\par{\bf Proof.}$\ \ $}

\def\includeE#1{{\lhook\kern-3.5pt\joinrel\smash{
    \mathop{\longrightarrow}\limits^{#1}}}}

\def\efor/{Example~\ref{E4}}

\def\ep{\par\bigskip}

\def\BL/{Baldwin--Lachlan}
\def\Bu/{Buechler}
\def\Hr/{Hrushovski}
\def\lm/{locally modular}
\def\wm/{weakly minimal}
\def\nm/{non--modular}
\def\tt/{totally transcendental}
\def\ss/{superstable}
\def\ud/{unidimensional}
\def\sm/{strongly minimal}

\def\abar{\overline{a}}

\def\bbar{\overline{b}}
\def\cbar{\overline{c}}
\def\dbar{\overline{d}}

\def\hbar{\overline{h}}

\def\xbar{\overline{x}}
\def\ybar{\overline{y}}

\def\tp{{\rm tp}}

\def\tr/{trivial}
\def\nt/{non--trivial}
\def\st/{strong type}

\def\TV/{Tarski--Vaught}

\def\sc/{sound construction}
\def\ac/{atomic construction}

\def\fal/{functional}
\def\upl/{unique parallel lines}
\def\chp/{categorical in a higher power}

\def\text#1{\ \hbox{#1}\ }
\def\subset{\subseteq}

\def\contains{\supseteq}
\def\<{\langle}
\def\>{\rangle}

\def\forces{\mathrel{\raise.4ex\hbox{$\scriptstyle \vert$\hskip-.5ex}\vdash}}

\def\K/{${\cal K}$}
\def\KM{{\cal K}}

\def\r{\restriction}
\def\rl{\r L}
\def\KMOD/{$\KM={\rm Mod}(T_1)\rl$}
\def\LK/{$L(\KM)$}
\def\S{{\rm Spec}(\KM)}
\def\al{\alpha}

\def\e{\emptyset}
\def\g{\gamma}

\def\l{\lambda}
\def\o{\omega}
\def\th{\theta}

\def\SN/{\hbox{$\S\ne\e$}}

\def\a1{\aleph_1}

\def\PCD/{${\rm PC}_{\a0}$}
\def\PC/{${\rm PC}_{\Delta}$}
\def\ra{\rightarrow}

\def\bg{\beth_\g}

\def\ss{\smallsetminus}

\def\SG/{\hbox{$\S\cap\bg$}}
\def\SNG/{\hbox{$\S\cap\bg\not =\e$}}

\def\A{{\cal A}}
\def\B{{\cal B}}
\def\C{{\cal C}}
\def\D{{\cal D}}

\def\P{{\cal P}}

\newtheorem{theorem}{Theorem}[section]

\newtheorem{lemma}[theorem]{Lemma}

\newremark{defn}[theorem]{Definition}
\newremark{example}[theorem]{Example}
\newremark{remark}[theorem]{Remark}
\author{
M. C. Laskowski\thanks{Partially supported by an NSF 
Postdoctoral Fellowship.}
\\Department of Mathematics\\
University of Maryland
\and
S. Shelah\\
Department of Mathematics\\
Hebrew University of Jerusalem
\and Department of Mathematics\\
Rutgers University
\thanks{The authors thank the U.S. Israel
Binational Science Foundation for its support of this project.
This is item 489 in Shelah's bibliography.}}
\date{December 14, 1992}
\title{On the existence of atomic models}

\begin{document}

\maketitle 
\begin{abstract}
We give an example of a countable theory $T$ such that for every
cardinal $\l\ge\aleph_2$ there is a fully indiscernible set $A$
of power $\l$ such that the principal types are dense over $A$,
yet there is no atomic model of $T$ over $A$.
In particular, $T(A)$ is a theory of size $\l$ where the principal
types are dense, yet $T(A)$ has no atomic model.
\end{abstract}
If a complete theory $T$ has an atomic model then the principal 
types are dense in the Stone space $S_n(T)$ for each $n\in\o$.
In \cite[Theorem 1.3]{Knight}, \cite[page 168]{Kueker} and 
\cite[IV 5.5]{Sh},  Knight, Kueker and Shelah independently showed that
the converse holds, provided that the cardinality of the
underlying language has size at most $\aleph_1$.

In this paper we build an example that demonstrates that
the condition  on the cardinality of the language is necessary.
Specifically, we construct a complete theory $T$ in a countable language
having a distinguished predicate $V$ such that if $A$ is any subset
of $V^M$ for any model $M$ of $T$, then the principal types are dense in
$S_n(T(A))$ for each $n\in\o$.  ($T(A)={\rm Th}(M,a)_{a\in A}$).
However, 
$T(A)$ has an atomic model if and only if $|A|\le\aleph_1$.

In fact, by modifying the construction (Theorem~\ref{strong}) we
may insist that there is a particular non-principal, complete type $p$ 
such that, for any subset $A$ of $V^M$,
$p$ is realized in every model of $T(A)$ if and only if $|A|\ge\aleph_2$.
With Theorem~\ref{general}
we show that the constructions can be generalized to larger cardinals.

\medskip

We first build a countable, atomic model in a countable language
having an infinite, definable subset of (total) indiscernibles.
Let $L$ be the language with unary predicate symbols
$U$ and $V$, a unary function symbol $p$,
and countable collections of 
binary function symbols $f_n$ and binary relation symbols $R_n$
for each $n\in\o$.
By an abuse of notation, $p$ and each of the $f_n$'s
will actually be partial functions.

For a subset $X$ of an $L$-structure $M$, define the closure of
$X$ in $M$, $cl(X)$, to be the transitive closure of
$cl_0(X)=\{f_n(b,c):b,c\in X, n\in\o\}$.  So $cl(X)$ is a subset of
the smallest substructure of $M$ containing $X$.

Let $K$ be the set of all finite $L$-structures $\A$ satisfying the
following eight constraints:
\def\labelenumi{\roman{enumi})}
\begin{enumerate}
\item $U$ and $V$ are disjoint sets whose union is the universe $A$;
\item $p:U\ra V$;
\item each $f_n:U\times U\ra U$;
\item for each $n$, $R_n(x,y)\ra (U(x)\wedge U(y))$;
\item the family $\{R_n:n\in\o\}$ partitions all of $U^2$ into
disjoint pieces;
\item for each $n$ and $m\ge n$, $R_n(x,y)\ra f_m(x,y)=x$;
\item if $x',y'\in cl(\{x,y\})$ and $R_n(x,y)$, then $\bigvee_{j\le n} R_j(x',y')$;
\item there is no $cl$-independent subset of $U$ of size 3
(i.e., for all \hbox{$x_0,x_1,x_2\in U$,} there is a permutation $\sigma$ of $\{0,1,2\}$
such that $x_{\sigma(0)}\in cl(\{x_{\sigma(1)},x_{\sigma(2)}\})$.)
\end{enumerate}

It is routine to check that $K$ is closed under substructures
and isomorphism and that $K$ contains only countably many isomorphism types.
We claim that $K$ satisfies the joint embedding property and the
amalgamation property.  As the proofs are similar, we only verify 
amalgamation.  Let $\A,\B,\C\in K$ with $\A\subset \B$, $\A\subset \C$
and $A=B\cap C$.  It suffices to find an element $\D$ of $K$ with
universe $B\cup C$ such that $\B\subset \D$ and $\C\subset \D$.
Let $\{b_0,\dots,b_{l-1}\}$ enumerate $U^{\B}$,
let $\{c_0,\dots,c_{m-1}\}$ enumerate $U^{\C}$ and let $k>l,m$ be large
enough that $\bigcup\{R_j^{\B}:j<k\}= (U^{\B})^2$
and
$\bigcup\{R_j^{\C}:j<k\}=(U^{\C})^2$.
For each $b\in B\ss A$ and $c\in C\ss A$, let $f_j(b,c)=c_j$ for 
$j< m$, and $f_j(b,c)=b$ for $j\ge m$  and let $R_k(b,c)$.  Similarly, for $j<l$ let
$f_j(c,b)=b_j$ and $f_j(c,b)=c$ for all $j\ge l$ and $R_k(c,b)$.
It is easy to check that $\D\in K$.

It follows (see e.g., \cite[Theorem 1.5]{KL}) that there is a countable, $K$-generic
$L$-structure $\B$.  That is, (*) $\B$ is the union of an increasing chain of
elements of $K$, (**) every element of $K$ isomorphically embeds into $\B$
and (***) if $j:\A\ra\A'$  is an isomorphism between 
finite substructures of $\B$ then there
is an automorphism $\sigma$ of $\B$ extending $j$.
Such structures are also referred to as homogeneous-universal structures.
Let $T$ be the theory of $\B$.  

We record the following facts about $\B$ and $T$:
First, $V^{\B}$ is infinite as there are elements $\A$ of $K$ with $V^{\A}$ arbitrarily
large;
$V^{\B}$ is a set of indiscernibles because of property  (***) and the fact that
any two $n$-tuples of distinct elements from $V^{\B}$
are universes of isomorphic substructures of $\B$; 
As every finite subset of $\B$ is contained in an element of $K$, it follows that
$cl(X)$ is finite for all finite subsets $X$ of $B$ and there is no
$cl$-independent subset of $U^{\B}$ of size 3;  Finally, $\B$ is atomic as
for any tuple $\bbar$ from $B$, property (***)
guarantees that the complete type of $\bbar$ is isolated
by finitely much of the atomic diagram of the smallest substructure $\A$ of $\B$
containing $\bbar$.  (If $n$ is least such that $\bigvee_{j\le n} R_j(a,a')$
for all $a,a'\in U^{\A}$ then we need only the reduct of the atomic diagram of
$\A$ to $L_n=\{U,V,p,R_j,f_j:j\le n\}$.)

\begin{lemma}  \label{0}
Let $\C$ be a model of $T$ and let $A$ be any subset of $V^{\C}$.
Then the principal types are dense over $A$.
\end{lemma}

\bp  Let $\th(\xbar,\abar)$ be any consistent formula, where $\abar$ is a
tuple of $k$ distinct elements from $A$.  Let $\bbar$ be any $k$-tuple of 
distinct elements from $V^{\B}$.  As the elements from $V$ are indiscernible,
$\B\models \exists\xbar\th(\xbar,\bbar)$.  Let $\cbar$ from $\B$ realize 
$\th(\xbar,\bbar)$.  Since $\B$ is atomic, there is a principal formula
$\phi(\xbar,\ybar)$ isolating $\tp(\cbar,\bbar)$.  It follows from indiscernibility
that $\phi(\xbar,\abar)$ is a principal formula such that
$\C\models \forall \xbar(\phi(\xbar,\abar)\ra\th(\xbar,\abar))$.
\ep

\begin{lemma}   \label{1}
Let $\C$ be an arbitrary model of $T$ and let $A\subset V^{\C}$.
If $\C$ is atomic over $A$  then:
\begin{enumerate}
\item $|U^{\C}|\ge |A|$;
\item $cl(X)$ is finite for all finite $X\subset U^{\C}$;
\item there is no $cl$-independent subset of size 3 in $U^{\C}$.
\end{enumerate}
\end{lemma}
\bp  

(i) $|U^{\C}|\ge |A|$ since for each $a\in A$, $p^{-1}(a)$ is non-empty.

(ii) As $\C$ is atomic, $\tp(X/A)$ is isolated by some formula
$\th(\xbar,\cbar)$, where $\cbar$ is a $k$-tuple of distinct elements from $A$.
As $V$ is indiscernible, $\th(\xbar,\bbar)$ is principal for any $k$-tuple
$\bbar$ of distinct elements from $\B$.  Choose $\dbar$ from $\B$ realizing
$\th(\xbar,\bbar)$ and suppose that $|cl(\dbar)|=l<\o$.
Then as $\th(\xbar,\abar)$ is principal, $\th(\xbar,\abar)$ implies 
$|cl(\xbar)|\le l$.

(iii) Assume $c_0,c_1,c_2\in U^{\C}$ are $cl$-independent.
As $\C$ is atomic over $A$, $\tp(\cbar)$ is principal, so let
$\th(\xbar,\abar)$ isolate $\tp(\cbar)$.  Again choose $\bbar$ from $V^{\B}$
and $\dbar$ from $U^{\B}$ such that  $\B\models \th(\dbar,\bbar)$.
But then $\dbar$ is a $cl$-independent subset of $U^{\B}$, which is a
contradiction.
\ep

We next record a well-known combinatorial lemma (see e.g., \cite{Erdos}).
An abstract closure relation on a set $X$ is a function $cl:\P(X)\ra\P(X)$
such that, for all subsets $A,B$ of $X$ and all $b\in X$,
$A\subset cl(A)$, $cl(cl(A))=cl(A)$, $A\subset B$ implies $cl(A)\subset cl(B)$ and
$b\in cl(A)$ implies there is a finite subset $A_0$ of $A$ such that
$b\in cl(A_0)$.

\begin{lemma}  \label{2}
For all ordinals $\al$ and all $n\in\o$, if $|X|\ge\aleph_{\al+n}$
and $cl$ is a closure relation on $X$ such that $|cl(A)|<\aleph_{\al}$
for all finite subsets $A$ of $X$.  Then $X$ contains a $cl$-independent
subset of size $n+1$.
\end{lemma}

\bp
Fix an ordinal $\al$.  We prove the lemma by induction on $n$.
For $n=0$ this is trivial, so assume the lemma holds for $n$.
Suppose $X$ has size at least $\aleph_{\al+n+1}$.  As $cl$ is finitely based
we can find a subset $Y$ of $X$, $|Y|=\aleph_{\al+n}$, such that
$cl(A)\subset Y$ for all $A\subset Y$.  Choose $b\in X\ss Y$.
Define a closure relation $cl'$ on $Y$ by
$cl'(A)=cl(A\cup\{b\})\cap Y$.
By induction there is a $cl'$-independent subset $B$ of $Y$ of size $n$.
It follows that $B\cup\{b\}$ is the desired $cl$-independent subset of $X$.
\ep

Note that by taking $\al=0$ in the lemma above, if $cl$ is a locally finite
closure relation on a set $X$  of size $\aleph_2$, then 
$X$ contains a $cl$-independent subset of size 3.

\begin{theorem}  \label{basic}
Let $A$ be a subset of $V^{\C}$ for an arbitrary model $\C$ of $T$.
Then $A$ is a set of
indiscernibles and the principal types over $A$ are dense, but there is an atomic
model over $A$ if and only if $|A|\le \aleph_1$.
\end{theorem}

\bp
The principal types are dense over $A$ by Lemma \ref{0}.  If $|A|\le\aleph_1$
then there is an atomic model over $A$ by e.g., Theorem~1.3 of \cite{Knight}.
However, if $|A|\ge\aleph_2$ then there cannot be an 
atomic model over $A$ by Lemmas~\ref{1}~and~\ref{2}.
\ep
Our next goal is to modify the construction given above so that there 
is a non-principal complete type $p$ that is realized in any model containing
$A$, provided that $|A|\ge\aleph_2$.
To do this, note that any atomic model of $T(A)$ is locally finite and
omits the type of a pair of elements from $U$ with every $R_n$ failing.
We shall enrich the language so as to code each of these by a single 1-type.

Let $L'=L\cup\{W,g,h\}\cup\{c_n:n\in\o\}$, where $W$ is a unary predicate, $g$ and $h$
are
respectively binary and ternary function symbols,
 and the $c_n$'s are new constant symbols.
Let $\B'$ be the $L'$-structure with universe $B\cup D$, where $D=\{d_n:n\in\o\}$
is disjoint from $B$, $W$ is interpreted as $D$, each $c_n$ is interpreted
as $d_n$, $g:U\times U\ra W$ is given by $g(a,b)=d_n$, where
$R_n(a,b)$ holds and $h(a,b,c)=d_n$ if and only if 
$n=|cl(\{a,b,c\})|$.
Let $T'$ be the theory of $\B'$.

Note that any automorphism $\sigma$ of $\B$ extends to an automorphism
$\sigma'$ of $\B'$, where $\sigma'\r D=id$.
It follows that $\B'$ is atomic and $V^{\B'}$ is an indiscernible set.
Let $p(x)$ be the non-principal type $\{W(x)\}\cup\{x\neq c_n:n\in\o\}$.
We claim that $p$ is complete.
This follows from the fact that for any $L'$-formula $\th(x)$, if
$n$ is greater than the number of terms occurring in $\th$ and
if $\th(x)$ is an $L'_n=\{U,V,W,g,h,p,R_l,f_l,
c_l:l<n\}$-formula then $\B'\models \th(c_i)\leftrightarrow \th(c_j)$
for all $i,j\ge n$.  (This fact can be verified by finding a back-and-forth
system 
${\cal S}=\{\langle \abar,\bbar\rangle:|\abar|=|\bbar|<n\}$
such that, for every $\langle\abar,\bbar\rangle\in {\cal S}$ and every
atomic 
$L'_n\cup\{R_i,f_i,c_i\}$-formula $\phi(\xbar)$, $|\xbar|=|\abar|$,
$$\B'\models \phi(\abar)\leftrightarrow\phi'(\bbar),$$
where $\phi'(\xbar)$ is the 
atomic
$L'_n\cup\{R_j,f_j,c_j\}$-formula generated from $\phi(\xbar)$ by replacing
each occurrence of $R_i,f_i,c_i$ by $R_j,f_j,c_j$, respectively.)

\begin{theorem}   \label{strong}
Let $A$ be a subset of $V^{\C'}$ for any model $\C'$ of $T'$.
Then $A$ is a set of indiscernibles and the principal types over $A$
are dense.  Further, if $|A|\le \aleph_1$ then there is an atomic model
over $A$, while if $|A|>\aleph_1$ then any model of $T'(A)$
realizes the complete type $p$.
\end{theorem}

\bp
The first two statements follow from the atomicity of $\B'$ and
the indiscernibility of $V^{\B'}$.  If $|A|\le\aleph_1$ then the
existence of the atomic model over $A$ follows from Knight's theorem.
So suppose $|A|\ge\aleph_2$ and let $\D'$ be any model of $T'$
containing $A$.  By examining the proof of
Lemma~\ref{2} it follows that either there are $a,b,c\in U^{\D'}$ such that
$cl(\{a,b,c\})$ is infinite or that there are $cl$-independent 
elements $a_0,a_1,a_2\in U^{\D'}$
(i.e., if $cl$ is a closure relation on $X$ and $|X|\ge\aleph_2$ 
then either $cl(x,y,z)$ is infinite for some $x,y,z\in X$ or there is an
independent subset of $X$ of size 3).

In the first case $h(a,b,c)$ realizes $p$. 
Now assume that the closure of any triple from $U^{\D'}$ is finite.
If, in addition, for every two elements $a,b$ from $U^{\D'}$
there were an integer $n$ such that $\D'\models R_n(a,b)$,
then $T'$ would ensure that there would not be any 3-element
$cl$-independent subset of $U^{\D'}$. (Under these assumptions there would be
only finitely many possibilities for the diagram of a triple
under the functions $\{f_i:i\in\o\}$ and no independent triple exists in $\B'$.)
Consequently, there must be a pair of elements $a,b$
from $U^{\D'}$ such that 
$\D'\models \neg R_n(a,b)$ for all $n\in\o$, so $g(a,b)$
realizes $p$.
\ep
We close with the following theorem demonstrating that the behavior between
$\aleph_1$ and $\aleph_2$ holds more generally between $\aleph_k$ and 
$\aleph_{k+1}$ for all $k\ge 1$.  The theorem is stated in its most basic
form to aid readability. 
We leave it to the reader to verify
that the strengthenings given in Theorems~\ref{basic} and \ref{strong}
(i.e., no atomic model over a given set or the non-atomicity being witnessed
by a specific complete type)
can be made to hold as well.

\begin{theorem}   \label{general}
For every $k$, $1\le k<\o$ there is a
countable theory $T_k$
such that $T_k$ has an atomic model of size $\aleph_\al$ if and only if $\al\le k$.
\end{theorem}

\bp
Fix $k$. Let $L_k=\{f_n,R_n:n\in\o\}$, 
where each $f_n$ is a  $(k+1)$-ary function and
each $R_n$ is a $(k+1)$-ary relation.  
Define $cl(X)$ to be the transitive closure of $cl_0(X)=\{f_n(a_0,\dots,a_k):a_i\in X\}$
and let $K$ be the set of all finite $L_k$-structures
satisfying the following constraints: 
\begin{enumerate}
\item $\{R_n:n\in\o\}$ partitions the $(k+1)$-tuples into
disjoint pieces;
\item for each $n$ and $m\ge n$, $R_n(x_0,\dots,x_k)\ra f_m(x_0,\dots,x_k)=x_0$;
\item if $x_0',\dots,x_k'\in cl(\{x_0,\dots,x_k\})$ and 
$R_n(x_0,\dots,x_k)$, then $\bigvee_{j\le n} R_j(x_0,\dots,x_k)$;
\item there is no $cl$-independent subset of size $k+2$.
\end{enumerate}

As before, there is a countable, $K$-generic $L_k$-structure $\B$.
Let $T_k=Th(\B)$.
Just as before, $\B$ is atomic, $cl$ is locally finite
on $\B$, $\bigvee_{n\in\o} R_n(\bbar)$ holds for all $(k+1)$-tuples $\bbar$ from $\B$
and there is no $(k+2)$-element $cl$-independent subset of $\B$.
Thus, the proof that there is no atomic model of $T_k$ of power $\l>\aleph_k$
is exactly analogous to the proof of Theorem~\ref{basic}.
What remains is to prove that there is an atomic model
of size $\aleph_k$.  
To help us, we quote the following combinatorial fact, which is a sort of converse
to Lemma~\ref{2}.
\begin{lemma}    \label{converse}
For every $k$, $1\le k<\o$ and every set $A$, $|A|\le \aleph_{k-1}$,
there is a family of functions
$\{g_n:A^k\ra A:n\in\o\}$ such that, letting $cl$ denote the transitive closure
under the $g_n$'s:
\begin{enumerate}
\item  $cl$ is locally finite;
\item  for all $\abar\in A^k$ there is an $n$ such that $g_m(\abar)\in\abar$ for all 
$m\ge n$;
\item there is no $cl$-independent subset of $A$ of size $k+1$.
\end{enumerate}
\end{lemma}

\bp
We prove this by induction on $k$.  If $k=1$, let $\{a_i:i<\al\le\aleph_0\}$ 
enumerate $A$.  Define $g_n$ by
$$g_n(a_i)=\cases{a_n&if $n <i$;\cr a_i&otherwise.\cr}$$

Now assume the lemma holds for $k$.  Let $\{a_i:i<\al\le\aleph_k\}$
enumerate $A$.  Define $g_n:A^{k+1}\ra A$ as follows:
Given $\abar=\langle a_{i_0},\dots,a_{i_k}\rangle\in A^{k+1}$,
let $i^*=\max\{i_0,\dots,i_k\}$.  
If $i^*=i_l$ for a unique $l<k+1$ then we can apply the inductive hypothesis to
the set $A_{i^*}=\{a_j:j<i^*\}$ and obtain a family of functions
$h_n:A_{i^*}^k\ra A_{i^*}$ satisfying the conditions of the lemma.
Now define $g_n(\abar)=h_n(\bbar)$, where $\bbar$ is the subsequence
of $\abar$ of length $k$ obtained by deleting $a_{i_l}$ from $\abar$.
On the other hand, if
there are $j<l<{k+1}$ such that
$i_j=i_l=i^*$, then simply let $g_n(\abar)=a_{i^*}$ for all $n\in\o$.
\ep

To show that there is an atomic model of $T_k$ of size $\aleph_k$,
as the principal formulas are dense and are $\Sigma_1$ 
it suffices
to show the following:
\medskip\par\indent ($\#$) If $\A$ is an $L_k$-structure of size at most
$\aleph_{k-1}$ such that every finitely generated substructure of $\A$
is an element of $K$ and $\phi(x,\dbar)$ ($\dbar$ from $A$)
is a principal formula consistent with $T_k$, then there is an extension
$\C\contains\A$ containing a witness to $\phi(x,\dbar)$ such that
every finitely generated substructure of $\C$ is in $K$.
\medskip
\par\indent
So choose $\A$ and $\phi(x,\dbar)$ as above.  We may assume that 
$cl(\dbar)=\dbar$.
We shall produce an
extension $\C\contains \A$ such that $C\ss A$ is finite,
$\C\models\exists x\phi(x,\dbar)$ and
every finitely generated substructure of $\C$ is in $K$.

Since $\phi(x,\dbar)$ is consistent with $T_k$ there is an element $\D$ of K
such that $\dbar$ embeds isomorphically into $\D$ and $\D\models \exists x\phi(x,\dbar)$.
Let $\{b_0,\dots,b_{l-1}\}$ enumerate the elements of $D\ss\dbar$.
We must extend the definitions of $\{f_n:n\in\o\}$ and $\{R_n:n\in\o\}$
given in $\A$ and $\D$
to $(k+1)$-tuples $\abar$ from $A\cup\{b_0,\dots,b_{l-1}\}$
so that every finitely generated substructure is in $K$.

We perform this extension by induction on $i<l$.  So fix $i<l$
and assume that $\{f_n,R_n:n\in\o\}$ have been extended to all $(k+1)$-tuples
from $A\cup\{b_j:j<i\}$ so that every finitely generated substructure is in $K$.
Let $\abar=\langle a_0,\dots,a_k\rangle$ be a $(k+1)$-tuple from $A\cup\{b_j:j\le i\}$
containing $b_i$
with at least one element not in $D$.  
If $b_i=a_s$ for some $s>0$, then let $f_n(\abar)=a_0$ for all $n\in\o$ and let
$R_0(\abar)$ hold.

On the other hand, if $b_i\neq a_s$ for all $s>0$, then $a_0=b_i$, so let
$\abar'=\langle a_1,\dots,a_k\rangle$, apply Lemma~\ref{converse}
to $A$ and $k$,  and let
$$f_n(b,\abar')=\cases{g_n(\abar')&if $g_n(\abar')\not\in\abar'$;\cr
                       b&otherwise.}$$
It is easy to verify that 
$cl$ is locally finite on $A\cup\{b_j:j\le i\}$ and that there is no
$cl$-independent subset of size $k+2$.  It is also routine to extend the
partition given by the $R_n$'s so as to preserve ii) and iii) in the definition
of $K$.  

Thus, we have succeeded in showing ($\#$), which completes the proof of 
Theorem~\ref{general}.
\ep


\begin{thebibliography}{99}
\bibitem{Erdos} P.\ Erd\"{o}s, A.\ Hajnal, A.\ Mate, P.\ Rado,
{\it Combinatorial Set Theory\/}, North Holland, Amsterdam, 1984.

\bibitem{Knight} J. Knight, Prime and atomic models,
{\it Journal of Symbolic Logic\/} {\bf 43} (1978) 385-393.

\bibitem{Kueker} D. W. Kueker, Uniform theorems in infinitary logic,
{\it Logic Colloquium~'77}, A. Macintyre, L. Pacholski, J. Paris (eds),
North Holland, 1978.

\bibitem{KL} D. W. Kueker and M. C. Laskowski, On generic structures,
{\it Notre Dame Journal of Formal Logic\/} {\bf 33} (1992) 175-183.

\bibitem{Sh} S. Shelah, {\it Classification Theory\/}, North Holland,
Amsterdam, 1978.

\end{thebibliography}
\end{document}